\documentclass[journal]{IEEEtran}

\usepackage{ifpdf}
 \ifpdf
 \else
 dvi code
 \fi
\usepackage{amsmath}
\usepackage{bm} % optional
\usepackage{lmodern}
\usepackage{moreverb}
\usepackage{color}
\usepackage{textcomp}
\usepackage{amsfonts,amssymb,amsbsy}
\usepackage{graphicx,empheq}
\usepackage{scrextend}
\ifCLASSINFOpdf

\else

\fi

\newcommand{\mc}{\mathcal}

\renewcommand{\theenumii}{\roman{enumii}}
\renewcommand{\theenumiii}{\alph{enumiii}}
\renewcommand{\thepage}{}

\def\e1{{\varepsilon_{11}}}
\def\b1{{\beta_{11}}}
\def\bp3{{\beta_{33}}}
\def\ep3{{\varepsilon_{33}}}
\newcommand{\mb}{\mathbf}

\def\mX{{\mathbb X}}
\def\Ltwo{{\mathbb L}^2 }

\newtheorem{thm}{Theorem}[section]

\newtheorem{prop}{Proposition}[section]

\newtheorem{Lem}{Lemma}[section]
\newtheorem{rmk}{Remark}[section]

\MHInternalSyntaxOn
\providecommand*\phantomword[3][c]{%
\mathchoice
{\MT_phantom_word:NNnn #1\displaystyle {#2}{#3}}%
{\MT_phantom_word:NNnn #1\textstyle {#2}{#3}}%
{\MT_phantom_word:NNnn #1\scriptstyle {#2}{#3}}%
{\MT_phantom_word:NNnn #1\scriptscriptstyle {#2}{#3}}%
}
\def\MT_phantom_word:NNnn #1#2#3#4{%
\@begin@tempboxa\hbox{$\m@th#2#4$}%
\setlength\@tempdima{\widthof{$\m@th#2#3$}}%
\hbox{\hb@xt@\@tempdima{\csname bm@#1\endcsname}}%
\@end@tempboxa}
\MHInternalSyntaxOff

\begin{document}
\title{Modeling and controlling an active constrained layer (ACL) beam  actuated by two voltage sources with/without magnetic effects}


\author{\IEEEauthorblockN{Ahmet \"Ozkan \"Ozer\IEEEauthorrefmark{1}\\}
\IEEEauthorblockA{\IEEEauthorrefmark{1}Department of Mathematics, Western Kentucky University, Bowling Green, KY 42101, USA}
%\thanks{As this paper was first submitted in September 2015, the author was a postdoctoral fellow and was affiliated by the Department of Mathematics at the University of Nevada-Reno.}
}

%\markboth{IEEE Transactions on Automatic Control}%
%{Shell \MakeLowercase{\textit{et al.}}: Bare Demo of IEEEtran.cls for IEEE Transactions on Magnetics Journals}



% for Transactions on Magnetics papers, we must declare the abstract and
% index terms PRIOR to the title within the \IEEEtitleabstractindextext
% IEEEtran command as these need to go into the title area created by
\IEEEtitleabstractindextext{%
\begin{abstract}
A fully dynamic three-layer active constrained layer (ACL) beam
% consisting of two piezoelectric elastic layers constraining a viscoelastic layer,
is modeled  for cantilevered boundary conditions by using a thorough variational approach. The Rao-Nakra thin compliant layer assumptions are adopted to model the sandwich structure, and all magnetic effects for the piezoelectric layers are retained. The piezoelectric layers are activated by two different voltage sources. When there are no ``mechanical" boundary forces acting in the longitudinal direction, it is shown that the system with certain parameter combinations is not uniformly strongly stabilizable by the $B^*-$type feedback controller, which is the total current accumulated at the electrodes for the piezoelectric layers. However, as the magnetic effects are ignored (electrostatic assumption), the closed-loop system with all mechanical feedback controllers is shown to be uniformly exponentially stable.
\end{abstract}

}

\maketitle
\IEEEdisplaynontitleabstractindextext
\IEEEpeerreviewmaketitle

\section{Introduction}
An active constrained layer (ACL) composite beam  consists of two piezoelectric layers and a constrained viscoelastic layer.  Each piezoelectric layer is actuated by  different voltage sources. When the electrodes of the piezoelectric layers are subjected to voltage sources, they  shrink or extend, and therefore, the whole ACL composite shrinks/extends or bends.
Accurately modeling the composite requires certain mechanical and electrical (and magnetic) assumptions for each layer. The middle layer is modeled by classical Mindlin-Timoshenko assumptions and the stiff piezoelectric layers are modeled by the Euler-bernoulli assumptions. Piezoelectric layers are traditionally modeled through the electrostatic assumption, and  all dynamic electrical effects and magnetic effects due to Maxwell's equations are ruled out, i.e. see  \cite{Smith}, and the references therein.  Since an ACL composite includes a piezoelectric layer, the corresponding models use the electrostatic assumption as well \cite{Trindade}. The reduced model (\cite{Baz,Lam,Shen}) is mostly either a Mead-Marcus-type \cite{Mead} or a Rao-Nakra-type model \cite{Rao}. For example, \cite{Baz} obtained a Mead-Marcus type model by neglecting the rotational inertia terms for the longitudinal dynamics and rotational inertia for the bending dynamics.  All of these models  reduce to the classical counterparts as in (\cite{Ditaranto,O-Hansen1,Mead}) once the piezoelectric strain is taken to be zero. On the other hand, the model obtained in \cite{Shen} through a consistent variational approach is more like a Rao-Nakra-type \cite{Rao}.

%\begin{figure}[h!tb]
%\centering
%\includegraphics[width=3.0in]{ACL-twopiezos.png}
%\caption{A three-layer ACL beam actuated by two different voltage sources $V^{\rm T}$ and $V_{\rm B}$.} \label{ACL}%\vspace{-0.3in}
 %\end{figure}
As the electrostatic assumption is adopted, it can be easily shown that a single piezoelectric beam model is exactly controllable and uniformly exponentially stabilizable for the $B^*-$type mechanical feedback, velocity of the beam at one end. The same type of phenomenon is observed for the ACL composite. For example in \cite{Baz}, it is shown that the time derivative of the  energy is nonnegative  as a mechanical damping is injected through the boundary of the piezoelectric layer.
 %It is not quite clear though if there are eigenvalues on the imaginary axis. As well, the exponential stabilizability of all eigen-modes are not proven to hold true.
 See other control strategies i.e. in \cite{Stanway}, and the references therein. However, even in the case of a single piezoelectric beam, as the dynamic effects are kept, a strongly coupled wave system is obtained for which it is shown that the model is not controllable and not uniformly exponentially stabilizable for almost  all combinations of material parameters with the $B^*-$type feedback, i.e.  the total current at the electrodes  \cite{O-M1}. In fact, there are no other feedback controllers to make system uniformly exponentially stable. It is worthwhile to mention that the closed-loop system obtained by supplying voltage to the electrodes and feeding back the total current is much easier and physical in terms of practical applications since measuring the total current at the electrodes is easier than measuring displacements or the velocity of the composite at one end of the beam, i.e. see (\cite{Chee,Miller}). In fact, with the same feedback controller, explicit polynomial decay estimates are obtained for more regular initial data \cite{Ozkan1}. These results for a single piezoelectric beam indicate that a similar type of controllability characteristics may be observed for composites involving piezoelectric layers, i.e. the closed-loop system for an ACL composite with a similar type of input-output mentioned above may be more physical, and moreover, the ACL composite beam models with the magnetic effects may lack stabilizability/controllability for some choices of their parameters.

To our knowledge, magnetic effects are ignored for any smart composite involving a piezoelectric layer. In this paper, a novel modeling strategy is proposed to obtain models of an ACL composite with/without the magnetic effects by using the Rao-Nakra assumptions \cite{Rao}. We further let the weight and the stiffness of the middle layer go to zero, i.e. see \cite{Hansen3}, to obtain reduced models which are the perturbations of the original models. Therefore, our models differ from the classical counterparts substantially not only due the inclusion of the magnetic effects but also due to the simplicity of the coupling between the longitudinal and bending dynamics, i.e. see (\cite{Baz,Lam}). A preliminary model has recently been obtained for an ACL composite with a single piezoelectric layer \cite{Ozkan2}. In the modeling process, more attention is paid to the piezoelectric layers  by using the full set of Maxwell's equations to incorporate the magnetic effects.  We prove the lack of uniform strong stabilization with the $B^*-$type feedback in the case where the mechanical boundary forces for the longitudinal dynamics are removed. Next, we consider the model without magnetic effects. This model is reduced to the \emph{``non-smart"} Rao-Nakra composite \cite{Hansen3}. To our knowledge, the boundary feedback stabilization or exact controllability for  the clamped-free Rao-Nakra beam model were never considered in the literature, i.e. see \cite{O-Hansen4} and the references therein, except the result for  the clamped-hinged boundary conditions where each equation for stretching and bending requires only one controller. In this paper, our results require two controllers for the bending equation. By adopting a similar argument of \cite{O-Hansen4}, it is shown that the reduced model is  a compact perturbation of the decoupled system consisting of a Rayleigh (Kirchhoff) beam equation \cite{B-Rao} and  wave equations. Therefore the coupled system is shown to be exponentially stable since the unique continuation result proves that there are no eigenvalues on the imaginary axis, see Lemma \ref{xyz}. Unlike the clamped-hinged ACL beam \cite{O-Hansen4}, our result requires using the four boundary conditions at one end which means that two controllers are essential to control the bending motion for the smart beam.
 %Our results use a similar technique  used for the same model with clamped-hinged boundary conditions \cite{O-Hansen3}.

\section{Modeling Preliminaries}
Consider an ACL beam occupying the region $\Omega=\Omega_{xy}\times (0, h)=[0,L]\times [-b,b] \times (0,h)$ at equilibrium where $\Omega_{xy}$ is a smooth bounded domain in the plane.
The total thickness $h$ is assumed to be small in comparison to the dimensions of $\Omega_{xy}$.
The beam  consists of two piezoelectric  layers and  a complaint layer.
 The layers are indexed from 1 to 3 from the bottom piezoelectric  layer to the top piezoelectric layer, respectively. Now we let $0=z_0<z_1<z_2<z_3=h, $ with $h_i=z_i-z_{i-1}, \quad i=1,2,3.$
We use the rectangular coordinates $(x,y)$ to denote points in $\Omega_{xy},$ and  $(X, z)$ to denote points in $\Omega = \Omega_{{\rm B}} \cup \Omega^{\rm ve} \cup \Omega^{{\rm T}} $, where $\Omega_{\rm B}=\Omega_{xy}\times (z_0,z_1), \Omega^{\rm ve}=\Omega_{xy}\times (z_1,z_2),$ and $\Omega^{\rm T}=\Omega_{xy}\times (z_2,z_3)$ are the reference configurations of the bottom piezoelectric, viscoelastic, and top piezoelectric layers, respectively.

For $(x,y,z)\in \Omega,$ let $U(x,y,z) = ({U_1,U_2,U_3})(x,y, z)$ denote the displacement vector of the point
(from reference configuration). In order to obtain a beam theory,  all displacements  are assumed to be independent of $y-$coordinate, and $U_2\equiv 0.$ The transverse displacements is  $w(x,y,z)= U_3(x)= w^i(x)$ for  any $i$ and $x\in [0,L].$ Define $u^i(x,y,z)=U_1(x,0,z_i)=u^i(x)$  for $i=0,1,2,3$ and  for all $x\in (0,L).$

We use the standard sandwich beam assumptions to model the ACL beams.  The modified constitutive equations for the piezoelectric layers are
\begin{eqnarray}
 \label{cons-eq50} &&\left\{
  \begin{array}{ll}
   T_{11}^{i}=\alpha^{i} S^{(i)}_{11}-\gamma^{i}\beta^{i} D_3^{i},\quad   E_1^{i}=\beta_{1}^{i}D_1^{i}  &\\
   E_3^{i}=-\gamma^i \beta^{i} S_{11}^{i}+\beta^{i} D_3^{i}, \quad i=1,3
  \end{array} \right.
\end{eqnarray}
where $T,S, D,E, c_{11}, \gamma, \epsilon,$ and $\beta$; are stress tensor, strain tensor, electrical displacement vector,
electric field intensity vector, elastic stiffness coefficients,
piezoelectric coefficients, permittivity coefficients, impermittivity
coefficients, and
 $\alpha^{i}=\alpha_1^{i} + (\gamma^{i})^2 \beta^{i}, ~~\alpha_1^{i}=c_{11}^{i},~~ \alpha^2=c_{11}^2$ and $\gamma^i=\gamma_{31}^i,~~\gamma_1^i=\gamma_{15}^i,~~\beta^i=\frac{1}{\varepsilon_{33}^i}, ~~\beta_{1}^i=\frac{1}{\varepsilon_{11}^i},$
and the middle layer is $  T_{11}=\alpha_1^2 S_{11},\quad T_{13}= 2G_{2} S_{13} $
where $G_2$ is the shear modulus of the viscoelastic layer, and refer to (\cite{Ozkan2,Smith}) for the description of piezoelectric and elasticity coefficients. Defining ${\hat z}_i = \frac{z^{i-1}+z^i}{2},$ the strain components for the viscoelastic layer, and the piezoelectric layers are respectively given by
\begin{eqnarray} \nonumber & S_{11}=\frac{\partial v^2}{\partial x}- (z-\hat z_i) \frac{\partial \psi^2}{\partial x}, \quad S_{13}=\frac{1}{2}\left(\psi^2+ w_x\right)=\frac{1}{2}\phi^2,&\\
 \label{strains} & S_{11}^i=\frac{\partial v^i}{\partial x}- (z-\hat z_i) \frac{\partial^2 w}{\partial x^2},\quad S^i_{13}=0,\quad i=1,3,&\\
 \label{defs1}   &\psi^i=\frac{u^i-u^{i-1}}{h_i}, \quad \phi^i= \psi^i + w_x, \quad v^i= \frac{u^{i-1}+u^i}{2}&
\end{eqnarray}
where $i=1,2,3,$ and in particular, $\phi^1=\phi^3=0, ~~\psi^1=\psi^3=-w_x, ~~\phi^2=\psi^2+ w_x.$ Here  $\psi^{i}$   can be viewed as the total rotation angles
 of the deformed filament within the $i^{\rm th}$ layer in the $x-z$ plane, and $\phi^i$ represent the shear angles within each layer, and $v^i$ represent the longitudinal displacements of the center line of the $i^{\rm th}$ layer. For details of the constitutive equations and parameters, the reader may refer to (\cite{Hansen3,O-M1}, or \cite{Ozkan2}).

\noindent{{\bf{Inclusion of  the electrical kinetic energy for the piezoelectric layers:}}}
% We follow the dynamic approach in \cite{O-M1} to use the full set of Maxwell's equations.
Let $\mathrm{B}^{\rm i}$ be the magnetic field $\textrm{B}^{i}(x)$ for the $i^{\rm th}$ piezoelectric layer for  $i=1,3,$ and have the only nonzero component $\textrm{B}_2^i(x).$ Assume also that the electric field of the $i^{\rm th}$ piezoelectric layer $E_1^i=0,$ and thus $D_1^i=0.$ Assuming that $D_3^i$ does not vary in the thickness direction $D_3^i(x,z,t)=D_3^i(x,t),$ it follows from the Amp\'{e}re-Maxwell equation that   $\textrm{B}_2^i=-\mu^i\int_0^x \dot D_3^i(\xi, z, t) ~d\xi$
where $\mu^i$ is the magnetic permeability of the $i^{\rm th}$ layer. Now we define $p^i:=\int_0^{x}  D_3^i(\xi,  t) ~d\xi$ to be the total electric charge at point $x.$ The magnetic energy for the $i^{\rm th}$ layer is
$ B^i =\frac{\mu^i}{2}\int_{\Omega} \left(\dot p^i\right)^2 dX.$
%\nonumber &=&  \frac{1}{2\mu^i}\int_{\Omega} (\textrm{B}_2^i)^2~dX\\

Assume that the beam is subject to a distribution of boundary forces $(\tilde g^1, \tilde g^3, \tilde g)$  along its edge $x=L,$ see \cite{Lagnese-Lions}.  Let $V^{\rm T}(t)$ and $V_{\rm B}(t)$ be the voltages applied at the electrodes of the piezoelectric layers, respectively. Then the total work done by all mechanical and electrical external forces is
\begin{eqnarray}
\nonumber\mb{W} &=& \int_0^L \left(- (p^1)_x V_{\rm B}- (p^3)_x V^{\rm T} \right)~dx +  g^1 v^1 (L) \\
\label{work-done} &&  + g^3 v^3 (L) +  g w(L)-M w_x(L).
\end{eqnarray}
where $M=m^1+m^3,$ $g^i=\int_{z_{i-1}}^{z_i} \tilde g^i ~dz,~ g=\int_0^h \tilde g ~dz,$ $m^i=\int_{z_{i-1}}^{z_i}(z-\hat z^i) \tilde g_i  ~dz$ for $i=1,3.$
The modified Lagrangian for the ACL beam is
$  \mb{L}= \int_0^T \left[\mb{K}-(\mb{P}+\mb{E})+\mb B +\mb{W}\right]~dt$
    where $\mb K=\sum_{i=1}^3 \mb K^i,$ $\mb P+\mb E=\mb P^2+\sum_{i=1,3}(\mb P^i + \mb E^i),$ and $\mb B=\mb B^1 + \mb B^3$ are the kinetic energy, the total stored energy, and the magnetic energy of the beam \cite{O-M1},
 \begin{eqnarray}
 \nonumber &&{\mb{K}}= \frac{1}{2} \int_0^L \left[\left(\sum_{i=1,3} \rho_ih_i(\dot v^i)^2\right) + \left(\sum_{i=1,3} \rho_ih_i \right) \dot w^2\right.\\
\nonumber   &&\quad\quad \left.+ \rho_2h_2 (\dot \psi^2)^2+ \left(\rho_1h_1 + \rho_3h_3\right) \dot w_x^2\right]~dx,\\
\nonumber  &&{\mb P+\mb E} = \frac{1}{2} \int_0^L \left[\alpha^2 h_2 \left((v^2_x)^2 + \frac{h_2^2}{12}(\psi_{x})^2 \right)\right.\\
\nonumber &&  + G_2h_2 \left(\phi^2\right)^2 + \sum_{i=1,3} \left(\alpha^i h_i  \left( (v_x^i)^2 +\frac{h_i^2}{12}(w_{xx})^2\right) \right. \\
\nonumber &&\left.\left.-2\gamma^i\beta^i h_i  v^i_x p^i_x + \beta^i h_i (p^i_x)^2\right)\right]~dx,\\
\label{k-energy}   && \mb{B}=\frac{1}{2}\int_{0}^L \sum_{i=1,3}\left(\mu^i h_i(\dot p^i)^2\right)~dx
\end{eqnarray}
where $\rho_i$ is the volume density of the $i-^{\rm th}$ layer.
Refer to (\cite{Hansen3,O-M1}) for the details.
\section{Rao-Nakra model and Hamilton's Principle }

By using (\ref{defs1}), $\{v^2,  \psi^2\}$ can be written as functions of $\{w,v^1,v^3\}.$ Thus,
we choose only $\{w,v^1,v^3\}$ as the state variables.  Let $H=\frac{h_1 + 2h_2+h_3}{2}.$ Application of Hamilton's principle by using cantilevered boundary conditions and  by setting the variation of admissible displacements $\{v^1,v^3, p^1, p^3, w\}$ of ${\mb L}$ to zero yields a highly coupled equations for bending and stretching of the whole composite. Thus, we  study the thin compliant layer Rao-Nakra model by letting $\rho_2, \alpha^2\to 0:$ %This
%approximation retains the potential energy of shear and transverse kinetic energy so that the model above reduces to
\begin{eqnarray}
 \label{dbas-mag} \left\{
  \begin{array}{ll}
 \rho_ih_i\ddot v^i-\alpha^i h_i v^i_{xx} +\gamma^i \beta^i h_i p^i_{xx}+ \kappa(i) G_2 \phi^2 = 0,   & \\
 \mu^i  h_i \ddot p^i   -\beta^i h_i  p^i_{xx} + \gamma^i \beta^i h_i v^i_{xx}= 0, \quad i=1,3, &\\
  m \ddot w - K_1 \ddot{w}_{xx} + K_2 w_{xxxx} - G_2 H \phi^2_x=0,&\\
 \phi^2=\frac{1}{h_2}\left(-v^1+v^3 + H w_x\right)&  \end{array} \right.
\end{eqnarray}
with the boundary and initial conditions for $i=1,3$
\begin{eqnarray}
 \label{d-son-mag} \left\{
  \begin{array}{ll}
 v^i(0)= p^i(0),   \alpha^i h_i v^i_{x}(L)-\gamma^i \beta^i h_i p^i_x(L)=g^i(t),&\\
 \beta^1 h_1 p^1_x(L) -\gamma^1\beta^1 h_1v^1_x(L)= -V_{\rm B}(t)& \\
 \beta^3 h_3 p^3_x(L) -\gamma ^3\beta^3 h_3v^3_x(L)= -V^{\rm T}(t)& \\
 w(0)=w_x(0)=0,\quad K_2 w_{xx}(L) = -M(t)&\\
 K_1 \ddot w_x(L) -K_2 w_{xxx}(L) + G_2 H \phi^2(L)=g(t)&\\
(v^1, v^3, p^1, p^3, w, \dot v^1, \dot v^3, \dot p^1, \dot p^3, \dot w)(x,0)\\
~~=(v^1_0,  v^3_0, p^1_0, p^3_0, w_0, v^1_1, v^3_1,  p^1_1, p^3_1, w_1)&
 \end{array} \right.
\end{eqnarray}
where $\kappa(i)={\rm sgn}(i-2), m=\sum_{i=1}^3\rho_ih_i, K_1=\frac{\rho_1 h_1^3}{12} +\frac{\rho_3 h_3^3}{12},$ and $K_2=\frac{\alpha^1 h_1^3}{12}+\frac{\alpha^3 h_3^3}{12}.$
Note that, different from a single piezoelectric model, the voltage controls $V^{\rm T}(t)$ and $V_{\rm B}(t)$ strongly couple the stretching and bending equations due the shear effect $\phi^2$ of the middle layer.

 Note that the role of the longitudinally applied mechanical boundary feedback controllers $g^1$ and $g^3$  are crucial to obtain a uniform exponential stabilization result (one control for each equation). However, once the mechanical boundary controllers for the stretching equations are removed, i.e. $g^1, g^3\equiv 0,$ the system is not even strongly stable by the $B^*-$type feedback for certain choices of material parameters. For $i=1,3$ define
 \begin{small}
\begin{eqnarray}
\nonumber \zeta_{\pm}^i = \frac{\sqrt{\frac{(\gamma^i)^2\mu^i }{\alpha^i_1}+\frac{\mu^i}{\beta^i}+\frac{\rho_i}{\alpha^i_1}\pm \sqrt{\left(\frac{(\gamma^i) ^2\mu^i}{\alpha^i_1 }+\frac{\mu^i}{\beta^i}
+\frac{\rho^i}{\alpha^i_1  }\right)^2-\frac{4\rho^i  \mu^i}{\beta^i \alpha^i_1 }}}}{\sqrt{2}} \\
\nonumber b_{\pm}^i = \frac{1}{2}\left(\gamma^i+\frac{\alpha^i_1}{\gamma^i\beta^i}-\frac{\rho_i}{\gamma^i \mu^i}\pm\sqrt{\left(\gamma^i+\frac{\alpha_i^1}{\gamma^i\beta^i}-\frac{\rho_i}{\gamma^i \mu^i}\right)^2+\frac{4\rho_i}{  \mu^i}} \right)
\end{eqnarray}
\end{small}
where $\zeta_+, \zeta_-,b_-,b_+\ne 0, b_-\ne b_+, \zeta_+\ne \zeta_-$ with $\zeta_+\zeta_-= \sqrt{\frac{\rho_i\mu^i}{\beta^i \alpha^i_1}}, b_- b_+=\frac{\rho_i}{\mu_i}.$
Let two piezoelectric beams have the same material properties, i.e. $\alpha^1_1=\alpha^3_1,\gamma^1=\gamma^3,$ etc. Under no mechanical longitudinal forces $g^1=g^3\equiv 0,$ we have the following result:
\begin{thm}
The system (\ref{dbas-mag})-(\ref{d-son-mag}) with $g^1,g^3\equiv 0$ is not strongly stable by the  feedback, i.e. $V^{\rm T}(t)=\frac{s_3\dot p^3(L)}{2h_3},$ $V_{\rm B}(t)=\frac{s_1\dot p^1(L)}{2h_1},$ $M(t)=-k_1\dot w_x(L)$ and $g(t)=k_2\dot w(L)$ for $s_1, s_2, k_1, k_2>0$ if $\frac{\zeta^i_{+}}{\zeta^i_{-}}=\frac{2n_i-1}{2m_i-1}$ for some $m_i,n_i\in\mathbb{N}, i=1,3.$
\end{thm}

\textbf{Proof:}  We prove that there are eigenvalues on the imaginary axis. Consider the eigenvalue problem corresponding to (\ref{dbas-mag})-(\ref{d-son-mag}) with $\lambda=\imath \tau$:
\begin{eqnarray}
 \label{dbas-eig-mag} &&\left\{
  \begin{array}{ll}
 \alpha^i_1 h_i v^i_{xx}-\gamma^i \beta^i h_i p^i_{xx} - \kappa(i) G_2 \phi^2 = -\tau^2 \rho_i h_i z^i,    & \\
    \beta^i h_i  p^i_{xx} - \gamma^i \beta^i h_i v^i_{xx}= -\tau^2 \mu^i  h_i  p^i, ~~i=1,3, &\\
    - K_2 w_{xxxx} + G_2 H \phi^2_x=-\tau^2 (m w  - K_1 w_{xx}),&
% \phi^2=\frac{1}{h_2}\left(-v^1+v^3 + H u_x\right)&
 \end{array} \right.
\end{eqnarray}
with the overdetermined boundary conditions
\begin{eqnarray}
 \nonumber   \left|w=w_x=v^i=p^i\right|_{x=0}= \left|v^i_x=p^i_x =p^i\right|_{x=L}=0,\\
  \label{d-son-eig-mag} w(L)=w_x(L)=w_{xx}(L)=w_{xxx}(L)=0, \quad i=1,3.
\end{eqnarray}
Let  $w(x)\equiv 0$ and
\begin{eqnarray}
\nonumber v^i(x)&=&k_1^i\frac{a^i_+ b^i_+\sin{(a^i_- x)} -a_-b^i_-\sin{(a^i_+ x)}}{a^i_+ a^i_-(b^i_+-b^i_-)}\\
\nonumber \quad  &&+k_2^i \frac{-a^i_+ \sin{(a^i_- x)} +a^i_-\sin{(a^i_+ x)}}{a^i_+ a^i_-(b^i_+-b^i_-)}, \\
\nonumber p^i(x)&=& k_1^i \frac{(a^i_+\sin{(a^i_- x)} -a^i_-\sin{(a^i_+ x)})b_1b_2}{a^i_+ a^i_-(b^i_+-b^i_-)}\\
\nonumber \quad   &&+ k_2^i\frac{a^i_-b^i_+ \sin{(a^i_+ x)} -a^i_+b^i_-\sin{(a^i_- x)}}{a^i_+ a^i_-(b^i_+-b^i_-)}
\end{eqnarray}
where   $a^i_+=\tau\zeta^i_+=\frac{(2n_i-1)\pi}{2L},$ $a^i_-=\tau\zeta^i_-=\frac{(2m_i-1)\pi}{2L}$ for some $m_i,n_i\in\mathbb{N},$ $i=1,3,$ $k_1^i=\frac{a^i_-b^i_+ \sin{(a^i_+ L)} -a^i_+b^i_-\sin{(a^i_- L)}}{a^i_+ a^i_-(b^i_+-b^i_-)},$ and $k_2^i=-\frac{(a^i_+\sin{(a^i_- L)} -a^i_-\sin{(a^i_+ L)})b^i_+b^i_-}{a^i_+ a^i_-(b^i_+-b^i_-)}. $
%\begin{eqnarray}\nonumber k_1^i=\frac{a^i_-b^i_+ \sin{(a^i_+ L)} -a^i_+b^i_-\sin{(a^i_- L)}}{a^i_+ a^i_-(b^i_+-b^i_-)},\\
%\nonumber  k_2^i=-\frac{(a^i_+\sin{(a^i_- L)} -a^i_-\sin{(a^i_+ L)})b^i_+b^i_-}{a^i_+ a^i_-(b^i_+-b^i_-)}.
%\end{eqnarray}
Here $v^1=v^3, p^1=p^3, w\equiv 0, \phi^2\equiv 0,$ and $(v^1,v^3,p^1,p^3,w)$ is the non-trivial solution of eigenvalue problem  (\ref{dbas-eig-mag})-(\ref{d-son-eig-mag}).
This implies that there are eigenvalues on the imaginary axis; $\left\{\pm\frac{\imath a^i_+}{\zeta^i_+}, \pm \frac{\imath a^i_-}{\zeta^i_-} \right\}.$ The conclusion  follows. $\square$

%Note that  the material parameters of the top and the bottom piezoelectric beams are the same; $\zeta_+^1=\zeta^3_+, \zeta_-^1=\zeta_-^3, b_+^1=b_+^3, b_-^1=b_-^3.$

\section{Stabilization without magnetic effects} First, assume that the magnetic energy for each layer is zero, i.e. $\mb B^i=0,$ and so,  $\ddot p^1=\ddot p^3\equiv 0$ in (\ref{dbas}). The electrostatic model is the well-known Rao-Nakra model in \cite{Hansen3}. The boundary stabilization problem is well studied in \cite{O-Hansen4} for the multi-layer beam clamped at the left end and hinged at the right end.
%This problem is relatively harder to handle in clamped-free boundary conditions case since the feedback operator is not a compact operator anymore.

  Finally, for $k_1,k_2,s_1,s_3> 0,$ analogous to \cite{Lagnese-Lions} and \cite{B-Rao},  we consider the following system
\begin{eqnarray}
 \label{dbas} &&\left\{
  \begin{array}{ll}
 \rho_ih_i\ddot v^i-\alpha^i_1 h_i v^i_{xx} + \kappa(i) G_2 \phi^2 = 0, \quad i=1,3,   & \\
   m \ddot w - K_1 \ddot{w}_{xx} + K_2 w_{xxxx} - G_2 H \phi^2_x=0,&\\
 \phi^2=\frac{1}{h_2}\left(-v^1+v^3 + H w_x\right)&  \end{array} \right.
\end{eqnarray}
with the boundary and initial conditions
\begin{eqnarray}
 \nonumber  &v^i(0)=w(0)=w_x(0)=0,~~K_2 w_{xx}(L) = -k_1 \dot w_x(L)&\\
   \nonumber &\alpha^i_1h_i v^i_x(L)=-s_i \gamma^i \dot v^i (L),\quad i=1,3,&\\
\nonumber & K_1 \ddot w_x(L) -K_2 w_{xxx}(L) + G_2 H \phi^2(L)=k_2\dot w(L)&\\
\label{d-son} &(v^1, v^3,  w, \dot v^1, \dot v^3,\dot w)(x,0)=(v^1_0,  v^3_0, w_0, v^1_1, v^3_1,   w_1).&
 \end{eqnarray}

\noindent {\bf{Semigroup well-posedness:}} \label{Sec-III}
Define  $ H^1_L(0,L)=\{\psi\in H^1(0,L): \psi(0)=0\},  H^2_L(0,L)=\{\psi\in H^2(0,L): \psi(0)=\psi_x(0)=0\},$
and  the complex linear spaces
\begin{eqnarray}
\nonumber \mX&=&\Ltwo(0,L), ~~ \mathrm V=\left(H^1_L(0,L)\right)^2 \times H^2_L(0,L),\\
\nonumber  \mathrm H&=& \mX^2 \times H^1_L(0,L), ~~ \mc{H} = \mathrm V \times \mathrm H.
\end{eqnarray}
so that $\mathrm V\subset \mathrm H \subset \mX^3\subset \mathrm H' \subset \mathrm V'.$
The  natural energy associated with (\ref{dbas})-(\ref{d-son})  is
\begin{eqnarray}
\nonumber && \mathrm{E}(t) =\frac{1}{2}\int_0^L \left\{\sum_{i=1,3}\left(\rho_i h_i  |\dot v^i|^2  +\alpha^i_1 h_i |v^i_x|^2\right) +m |\dot w|^2  \right.\\
\label{Energy-nat}  && \quad \left.     +K_1 |\dot w_x|^2 + K_2 |w_{xx}|^2+G_2 h_2 |\phi^2|^2 \right\}~ dx.~
\end{eqnarray}
 This motivates definition of the inner product on $\mc H$
{ \small{\begin{eqnarray}
\nonumber && \left<\left[ \begin{array}{l}
 u_1 \\
 \vdots \\
  u_6
 \end{array} \right], \left[ \begin{array}{l}
 v_1 \\
 \vdots \\
  v_6
 \end{array} \right]\right>_{\mc H}= \left<\left[ \begin{array}{l}
 u_4\\
 u_5\\
 u_6
 \end{array} \right], \left[ \begin{array}{l}
 v_4\\
 v_5\\
 v_6
 \end{array} \right]\right>_{\mathrm H}\\
 \nonumber && \quad\quad\quad+  \left<\left[ \begin{array}{l}
 u_1 \\
 u_2 \\
 u_3
 \end{array} \right], \left[ \begin{array}{l}
 v_1 \\
 v_2\\
 v_3
 \end{array} \right]\right>_{\mathrm V}\\
\nonumber && =\int_0^L \left\{\rho_1 h_1  u_4 {\dot {\bar v}}_4 + \rho_3 h_3  u_5 {\dot {\bar v}}_5 + m \dot u_6{\dot {\bar v}_6}+ K_1 (u_6)_x(\bar v_6)_x\right.  \\
\nonumber  &&+ \alpha^1_1 h_1  (u_1)_{x} (\bar v_1)_x + \alpha^3_1 h_3  (u_2)_{x} (\bar v_2)_x+  K_2 (u_3)_{xx} (\bar v_3)_{xx} \\
\nonumber &&  \left.+ \frac{G_2}{h_2}  (-u_1+u_2 + H(u_3)_x)(-\bar v_1+\bar v_2 + H(\bar v_3)_x) \right\}~dx.
 \end{eqnarray}}}
  Obviously, $\langle \, , \, \rangle_{\mc H} $ does indeed define an inner product, with the induced energy norm  since  the term $\|-u_1+u_2 + H(u_4)_x\|_{L^2(0,L)}$ is coercive, see \cite{Hansen3} for the details.

Let $\vec y=(v^1, v^3, w)$ be the smooth solution of the system of (\ref{dbas})-(\ref{d-son}).  Assuming the homogenous problem, all external forces are zero, multiplying the equations in (\ref{dbas}) by $\tilde y_1, \tilde y_3,\in H^1_L(0,L)$ and $ \tilde y \in H^2_L(0,L),$ respectively, and integrating by parts yields
\begin{eqnarray}
 \nonumber  &&\int_0^L \left(\rho_ih_i\ddot v^i \tilde y_i + \alpha^i_1 h_i v^i_{x} (\tilde y_i)_x + \kappa(i) G_2 \phi^2 \tilde y_i\right)dx \\
 \nonumber && \quad\quad\quad\quad + s_i  \dot v^i(L) \tilde y_i(L)= 0,   ~~~i=1,3, \\
  \nonumber && \int_0^L \left( m \ddot w \tilde y + K_1 \ddot{w}_{x} \tilde y_x   + K_2 w_{xx} \tilde y_{xx}- G_2 H \phi^2_x \tilde y_x \right)dx\\
\label{var} &&  +k_1 \dot w_x(L) \tilde y_x(L) + k_2 \dot w(L) \tilde y(L) =0.\quad\quad
\end{eqnarray}
Now define the linear operators
\begin{eqnarray}
\nonumber \left<Ay,\psi\right>_{\mathrm V'\times \mathrm V}=(y,\psi)_{\mathrm V \times \mathrm V},  \forall y,\psi\in \mathrm V\\
\nonumber \left<B_0\vec  y, \vec \psi\right>_{\mathrm H'\times \mathrm H}= \left[ \begin{array}{c}
 0_{2\times 1}\\
 k_2y_3(L) \psi_3(L)
 \end{array} \right], \forall \vec y,\vec \psi\in \mathrm H\\
\label{ops}\left<D_0\vec  y, \vec \psi\right>_{\mathrm H'\times \mathrm H}= \left[ \begin{array}{c}
 s_1y_1(L) \psi_1(L) \\
 s_3y_2(L) \psi_2(L) \\
 k_1(y_3)_x(L) (\psi_3)_x(L)
 \end{array} \right],  \forall \vec y,\vec \psi\in \mathrm V.
\end{eqnarray}
Let $\mc M : H^1_L (0,L) \to (H^1_L(0,L))'$ be a linear operator defined by
\begin{eqnarray}\label{def-M}\left< \mc M \psi, \tilde \psi \right>_{(H^1_L(0,L))', H^1_L(0,L)}= \int_0^L (m \psi \tilde {\bar \psi} + K_1 \psi_x \tilde {\bar\psi}_x) dx.
\end{eqnarray}
From the Lax-Milgram theorem $\mc M$ and $A$ are  canonical  isomorphisms from $H^1_L(0,L)$ onto  $(H^1_L(0,L))'$ and from $V$ onto $V',$ respectively.
Assume that $A y \in V',$ then we can formulate the variational equation above into the following form
$ M\ddot y  + A y +  D_0\dot y +  B_0\dot y=0$
where $M=\left[\rho_1 h_1 I ~~\rho_3 h_3 I~~ \mc M\right]$ is an isomorphism from $\mathrm H$ onto $\mathrm H'.$
Next we introduce the linear unbounded operator by
\begin{equation}
\label{defA}\mc A: {\text{Dom}}(A)\times V\subset \mc H \to \mc H
\end{equation}
where $\mc A= \left[ {\begin{array}{*{20}c}
   O_{3\times 3}  & -I_{3\times 3} \\
       M^{-1}A  &  M^{-1}  D_0  \\
\end{array}} \right]$
with ${\rm {Dom}}(\mc A) = \{(\vec z, \vec {\tilde z}) \in V\times V, A\vec z \in \mathrm V' \},$
and if ${\rm Dom}(\mc A)'$ is the dual of ${\rm Dom}(\mc A)$ pivoted with respect to $\mc H,$  we define the control operator $B$
 \begin{eqnarray}
\label{defb_0} \quad B \in \mathcal{L}(\mathbb{C} , {\rm Dom}(\mc A)'), ~ \text{with} ~ B=   \left[ \begin{array}{c} 0_{3\times 1} \\  M^{-1}B_0 \end{array} \right].
%\label{defb_0} \quad D \in \mathcal{L}(\mathbb{C} , {\rm Dom}(\mc A))', ~ \text{with} ~ D=   \left[ \begin{array}{c} 0_{3\times 1} \\ D_0 \end{array} \right]
 \end{eqnarray}

Writing $\varphi=[v^1, v^3, w, \dot v^1,  \dot v^3,  \dot w]^{\rm T},$ the control  system (\ref{dbas})-(\ref{d-son})  with the feedback controllers can be put into the  state-space form
\begin{eqnarray}
\label{Semigroup}
\dot \varphi + \mc A \varphi + B\varphi =0, \quad\varphi(x,0) =  \varphi ^0.
\end{eqnarray}

 \begin{Lem} \label{skew}The operator $\mc{A}$  defined by (\ref{defA}) is maximal monotone in the energy space $\mc H,$
 and ${\rm Range}(I+\mc A)=\mc H.$
\end{Lem}

\textbf{Proof:} Let $\vec z =\left[ \begin{array}{c}
\vec z_1 \\
\vec z_2 \end{array} \right]\in {\rm Dom}(\mc A).$ A simple calculation using integration by parts and the boundary conditions yields
{\small{\begin{eqnarray}
\nonumber &&\left<\mc A  \vec z, \vec z\right>_{\mc H\times \mc H} = \left<  \left[ \begin{array}{c}
-\vec z_2 \\
M^{-1} (A \vec z_1+D_0\vec z_2) \end{array} \right], \left[ \begin{array}{c}
\vec z_1 \\
\vec z_2 \end{array} \right]\right>_{\mc H}\\
\nonumber && \quad\quad= \left<-\vec z_2, \vec z_1\right>_{V\times V}  + \left<M^{-1}\left(A \vec z_1+D_0\vec z_2  \right), \vec z_2\right>_{H\times H}\\
\label{why1} &&\quad\quad = -\overline{\left<A \vec z_1, \vec z_2\right>_{V'\times V}}  + \left<A\vec z_1 + D_0\vec z_2, \vec z_2\right>_{H'\times H}.
\end{eqnarray}
}}
Since $\vec z \in {\rm Dom}(\mc A) $, then $A\vec z _1+ D_0\vec z_2\in V'$ and $\vec z_2 \in V$ so that
\begin{eqnarray}\nonumber\left<A\vec z_1 +  D_0\vec z_2, \vec z_2\right>_{H'\times H}=\left<A\vec z_1+ D_0\vec z_2, \vec z_2\right>_{V'\times V}\\
\label{why} \quad =\left<A\vec z_1, \vec z_2\right>_{V'\times V}+\left< D_0z_2, \vec z_2\right>_{V'\times V}.
\end{eqnarray}
Hence plugging (\ref{why}) in (\ref{why1}) yields
${\rm Re}\left<\mc A  \vec z, \vec z\right>_{\mc H\times \mc H}=\left< D_0\vec z_2, \vec z_2\right>_{V'\times V}\ge 0$ by (\ref{ops}).
We next verify the range condition. Let $\vec z=\left[ \begin{array}{c}
\vec z_1 \\
\vec z_2 \end{array} \right]\in \mc H.$ We prove that there exists a $\vec y=\left[ \begin{array}{c}
\vec y_1 \\
\vec y_2 \end{array} \right]\in \mc {\rm Dom} (\mc A )$ such that $(I+\mc A) \vec y=\vec z.$
A simple computation shows that proving this is equivalent to proving ${\rm Range} (M+A + D_0)=H',$ i.e., for every $\vec f \in H'$   there exists a unique solution $\vec z \in H$ such that $(M+A+D_0)\vec z=\vec f.$
This obviously follows from the Lax Milgram's theorem. $\square$

\begin{prop} The operator $B$ is a monotone compact operator on $\mathrm H.$
\end{prop}

\textbf{Proof:} Let $\left[ \begin{array}{c}
\vec y\\
\vec z \end{array} \right]\in \mathrm H.$ Then $\left<B \left[ \begin{array}{c}
\vec y\\
\vec z \end{array} \right], \left[ \begin{array}{c}
\vec y\\
\vec z \end{array} \right] \right>_{\mc H}=  k_2 |z_3(L)|^2.$
%\left<M^{-1}B_0 \vec z,\vec z\right>_{\mathrm H, \mathrm H}\\
%\nonumber  &=&  \left<B_0 \vec z,\vec z\right>_{\mathrm H', \mathrm H}\\
%\end{eqnarray*}
The compactness follows from the fact that $M^{-1}$ is a canonical isomorphism from $\mathrm H$ to $\mathrm H',$ and the fact that $B$ is a rank-one operator, hence compact from $\mathrm H$ to $\mathrm H'.$ $\square$

\subsection{Description of ${\rm Dom}(\mc A)$}

\begin{prop}\label{prop-dom}Let  $\vec u=(\vec y, \vec z)^{\rm T}\in \mc H.$ Then $\vec u \in {\rm Dom}(\mc A)$  if and only if the following conditions hold: $\vec y\in (H^2(0,L)\cap H^1_L(0,L))^2 \times (H^3(0,L) \cap H^2_L(0,L)),$
$\vec z\in V ~{\rm such ~ that}~(y_1)_{x}=(y_2)_{x}=(y_4)_{xx}\left. \right|_{x=L}=0.$
Moreover, the resolvent of $\mc A$ is compact in the energy space $\mc H.$
\end{prop}
\vspace{0.1in}

\textbf{Proof:} Let $\vec {\tilde u}= \left( \begin{array}{c}
\vec {\tilde y} \\
\vec {\tilde z}\\
 \end{array} \right)\in \mc H$ and $\vec {u}= \left( \begin{array}{c}
\vec {y} \\
\vec {z}\\
 \end{array} \right)\in {\rm Dom}(\mc A)$ such that $\mc A \vec u=\vec {\tilde u}.$ Then we have
 $ -\vec z= \vec {\tilde y} \in V, \quad   A\vec y + D_0 z=M\vec{\tilde z},$ and therefore,
\begin{eqnarray}\label{sal}\left< \vec y, \vec \varphi\right>_{V}= \left< \vec {\tilde z}, \vec \varphi \right>_{H} ~~{\rm for ~~all ~~} \vec \varphi \in V.\end{eqnarray}
Let  $\vec \psi=[\psi_1,\psi_2,\psi_3]^{\rm T} \in (C_0^\infty(0,L))^4.$  We define $\varphi_i=\psi_i$ for $i=1,2,$ and $\varphi_3=\int_0^x \psi_3 (s)ds.$ Since $\vec \varphi\in V,$ inserting $\vec \varphi$ into the above equation yields
\begin{eqnarray} \nonumber &&\int_0^L \left\{  -\alpha^1_1 h_1  (y_1)_{xx} \bar \psi_1  -\alpha^3_1 h_1  (y_2)_{xx} \bar \psi_2- K_2 (y_3)_{xxx} \bar \psi_3 \right.\\
\nonumber && \left. + \frac{G_2}{h_2}  (-y_1+y_2 + H(y_3)_x)(-\bar \psi_1+\bar \psi_2 + H(\bar \psi_3)_x) \right\}~dx\\
\nonumber && +s_1 (z_1)_x(L)(\psi_1)_x(L) + s_3 z_2(L)\psi_2(L) \\
  \nonumber &&+ k_1(z_3)_{x}(L) \psi_3(L)= \int_0^L \left\{ \left(\int_1^x m \tilde z_4  ds + K_1 (\tilde z_4)_x\right) \bar \psi_4 \right.\\
  \nonumber && \left.+\rho_1 h_1 \tilde z_1 \bar \psi_1+ \rho_3 h_3 \tilde z_2 \bar \psi_2 + \mu h_3 \tilde z_3 \bar \psi_3 \right\} ~dx
 \end{eqnarray}
for all $\vec \psi\in (C_0^\infty(0,L))^3.$
Therefore it follows that $\vec y\in (H^2(0,L)\cap H^1_L(0,L))^2 \times (H^3(0,L) \cap H^2_L(0,L)).$

Next let $\vec \psi \in \mathrm H.$  We define
\begin{eqnarray} \label{dumber}\varphi_i=\int_0^x \psi_i(s)ds,\quad i=1,\ldots,3.
\end{eqnarray}
Obviously $\vec \varphi \in \rm V. $ Then plugging (\ref{dumber}) into (\ref{sal}) yields
\begin{eqnarray}&& \nonumber  0=(\alpha^1_1 h_1  (y_1)_{x}(1) + s_1 z_1(L)) \bar \psi_1(1)+\alpha^3_1 h_3  (y_2)_{x}(1)\\
 \nonumber &&+ s_3 z_2(L) \bar \psi_2(1)+ (k_1 (y_3)_{xx}(L) + k_1 (z_3)_x(L)) (\bar \psi_3)_x(L)
\end{eqnarray}
for all $\psi\in \mathrm H.$ Hence,
\begin{eqnarray}\nonumber &&\alpha^1_1 h_1  (y_1)_{x}(1) + s_1 z_1(L)=\alpha^3_1 h_3  (y_2)_{x}(1) + s_3 z_2(L)=0,\\
\nonumber && k_1 (y_3)_{xx}(L) + k_1 (z_3)_x(L)=0.
\end{eqnarray}

Now let  $\vec y=\left[ \begin{array}{c}
\vec y_1 \\
\vec y_2 \end{array} \right]\in \mc {\rm Dom} (\mc A )$  and $\vec z=\left[ \begin{array}{c}
\vec z_1 \\
\vec z_2 \end{array} \right] $such that $(I+\mc A) \vec y=\vec z.$ By Proposition \ref{prop-dom} and  Lemma \ref{skew},  the compactness of the resolvent follows. $\square$

%%% birinci ve ikinic denklem bir kere turevle ...carpan olarak xz_xx al.....diger equation Rao gib.Sonra birinci ve ikinci denklemi z_xx ile carp.

\begin{Lem} \label{xyz} The eigenvalue problem
\begin{eqnarray}
 \label{dbas-eig} &&\left\{
  \begin{array}{ll}
 \alpha^i_1 h_i z^i_{xx} - \kappa(i) G_2 \phi^2 = \lambda^2 \rho_i h_i z^i, ~~i=1,3   & \\
    - K_2 u_{xxxx} + G_2 H \phi^2_x=\lambda^2 (m u  - K_1 u_{xx}),&
 %\phi^2=\frac{1}{h_2}\left(-z^1+z^3 + H u_x\right)&
 \end{array} \right.
\end{eqnarray}
with the overdetermined boundary conditions
\begin{eqnarray}
 \nonumber &&  u(0)=u_x(0)=z^i(0)=z^i(L)=z^i_x (L)=0, \quad i=1,3,\\
  \label{d-son-eig} && u(L)=u_x(L)=u_{xx}(L)=u_{xxx}(L)=0
\end{eqnarray}
has only the trivial solution.
\end{Lem}
\vspace{0.1in}
\textbf{Proof:} Now multiply the equations in (\ref{dbas-eig}) by $x\bar u_x-3\bar u,$ $x \bar z^1_x-2\bar z^1,$ and $ x\bar z^3_x-2\bar z^3,$ respectively, integrate by parts on $(0,L),$ and add them up:
\begin{eqnarray}
 \label{CH3-mult1-20}  \begin{array}{ll}
  \int_0^L\left\{ -\alpha^1_1 h_1|z^1_x|^2 -\alpha^3_1 h_3 |z^3_x|^2   -3\rho_1 h_1 \lambda^2|z^1|^2\right.\\
   -3\rho_3 h_3  \lambda^2 |z^3|^2 -4m \lambda^2 |u|^2 -  2 K_1 \lambda^2|u_x|^2   \\
   - G_2 h_2 \bar \phi^2 (z \phi^2_x)-3G_2 h_2|\phi^2|^2-K_2 \bar u_{xxxx} (x u_x) \\
     + \alpha^1_1 h_1 \bar z^1_{xx} (x z^1_x) + \alpha^3_1 h_3 \bar z^3_{xx}(x z^3_x) -\rho_1h_1 \lambda^2\bar z^1 (xz^1_x)\\
 \left. -\rho_3h_3 \lambda^2 \bar z^3 (x z^3_x)  -\lambda^2 (m\bar u-K_1 \bar u_{xx})(x u_x)\right)~dx=0.
 \end{array}
\end{eqnarray}
Now consider the conjugate eigenvalue problem corresponding to (\ref{dbas-eig})-(\ref{d-son-eig}). Now multiply the equations in the conjugate problem by $x u_x+2 u,$ $x z^1_x+3 z^1,$ and $ x z^3_x+3 z^3,$ respectively, integrate by parts on $(0,L),$ and add them up:
\begin{eqnarray}
\label{CH3-mult1-21}  \begin{array}{ll}
\int_0^L\left\{ 3 \left(\alpha^1_1 h_1 |z^1_x|^2 +\alpha^3_1 h_3|z^3_x|^2\right)   +3\bar\lambda^2\left(\rho_1 h_1 |z^1|^2\right.\right.\\
 \left.+\rho_3 h_3 |z^3|^2 \right) +2 \bar \lambda^2 \left(m  |u|^2 +   K_1 |u_x|^2\right) +2K_2 |u_{xx}|^2  \\
 + G_2 h_2 \bar \phi^2 (x     \phi^2_x)+3G_2 h_2|\phi^2|^2+K_2 \bar u_{xxxx} (x u_x) \\
  - \alpha^1_1 h_1 \bar z^1_{xx} (x z^1_x) - \alpha^3_1 h_3 \bar z^3_{xx}(x z^3_x) +\rho_1h_1 \bar \lambda^2\bar z^1 (xz^1_x)\\
  \left. +\rho_3h_3 \bar \lambda^2 \bar z^3 (x z^3_x)  +\bar \lambda^2 (m\bar u-K_1 \bar u_{xx})(x u_x)\right)~dx=0.
 \end{array}
\end{eqnarray}
Finally, adding (\ref{CH3-mult1-20}) and (\ref{CH3-mult1-21}),considering only the real part of the expression above and  all eigenvalues are
located on the imaginary axis, i.e. $\lambda=\mp \imath \nu,$ yields
\begin{eqnarray}
\nonumber \int_0^L\left( K_2 |u_{xx}|^2+ m \nu^2|u|^2+\sum_{i=1,3}(\alpha^i_1 h_i |z^i_x|^2)   \right)~dx=0.
 \end{eqnarray}
This implies that $u=z^1=z^3\equiv 0$ by  (\ref{d-son-eig}). In the case of $\lambda=0,$ we have
%\begin{eqnarray}
% \label{dbas-eig-zero} &&\left\{
 % \begin{array}{ll}
 %\alpha^i_1 h_i z^i_{xx} + \kappa(i) G_2 \phi^2 = 0, \quad i=1,3,  & \\
 %  - K_2 u_{xxxx} + G_2 H \phi^2_x=0.&
 %\phi^2=\frac{1}{h_2}\left(-z^1+z^3 + H u_x\right)&
 %\end{array} \right.
%\end{eqnarray}
%We simplify the above equation to
\begin{eqnarray}
 \label{dbas-eig-zero} &&\left\{
  \begin{array}{ll}
 -\phi^2_{xx} + \left(\frac{1}{\alpha^1_1 h_1} + \frac{1}{\alpha^3_1h_3}\right)\phi^2=-Hu_{xxx}\\
   - K_2 u_{xxxx} + G_2 H \phi^2_x=0,& \end{array} \right.
\end{eqnarray}
Since the operator $J=-D_x^2 + \left(\frac{1}{\alpha^1_1 h_1} + \frac{1}{\alpha^3_1h_3}\right)I$ is a non-negative operator on $H^2_*(0,L)=\{\psi\in H^2(0,L): \psi(0)=\psi_x(L)=0\},$ we obtain
$- K_2 u_{xxxx} -G_2 H^2 (J^{-1} u_{xxx})_x=0,$
and since $K_2D_x^4 + D_x J^{-1}D_x^3 $ is a positive operator on its domain, $u=\phi^2=0.$ And therefore, $u=z^1=z^3\equiv 0.$

\begin{thm} \label{strong}The semigroup generated by $(\mc A +B)$ is strongly stable in $\mc H.$
\end{thm}

\textbf{Proof:} We know that the system is dissipative, i.e. $\left<(\mc A + B) \left[ \begin{array}{c}
\vec z\\
\dot {\vec z} \end{array} \right], \left[ \begin{array}{c}
\vec z\\
\dot {\vec z} \end{array} \right]\right>_{\mathrm H}\le 0$. This result together with Lemma \ref{xyz} imply that there are no eigenvalues on the imaginary axis. The conclusion follows. $\square$

\begin{rmk} The method used to prove  Lemma \ref{xyz} is not valid once we remove the boundary condition $u(L)=0.$
%It is different from the result obtained in \cite{O-Hansen3} since any of the boundary conditions (clamped-clamped, clamped-hinged, and hinged-hinged) considered there involves  $u(L)=0.$
%The technique of multiplying the bending equation by $xu_{xxx}$ as in \cite{B-Rao} does not apply here due to the existence of the $\phi^2$ term in each equation. This may be done if one can show that the \emph{coupled} system has the Riesz basis property as in \cite{Guo,Wang}.
An analogous result obtained in \cite{O-Hansen4} was for either clamped, hinged, or mixed boundary conditions. It does require $u(L)=0.$ The interesting question is whether Theorem \ref{strong} without $u(L)=0$ is obtained for the case $k_1\ne 0, k_2\equiv 0.$  %This is an open question. If this result obtained, it is possible to exactly control/uniformly stabilize the system (\ref{dbas-mag})-(\ref{d-son-mag}) without the magnetic effects and with only three controllers.
%As a side comment, the multiplier method proposed in \cite{B-Rao} to reduce the number of controllers to one for the bending motion is not applicable for the coupled system.
\end{rmk}

Now we  consider the decomposition $\mc{A}+B=(\mc{A}_d+B)+ \mc {A}_{\phi}$ of the semigroup generator of the original problem (\ref{defA}) where $\mc{A}_d+B$ is the semigroup generator of the decoupled system, i.e. $\phi^2\equiv 0$ in (\ref{dbas})-(\ref{d-son}),
\begin{eqnarray}
 \label{dbas-dcpld} &&\left\{
  \begin{array}{ll}
 \rho_ih_i\ddot v^i-\alpha^i_1 h_i v^i_{xx} = 0,  \quad i=1,3, & \\
   m \ddot w - K_1 \ddot{w}_{xx} + K_2 w_{xxxx} =0,&
 %\phi^2=\frac{1}{h_2}\left(-v^1+v^3 + H w_x\right)&
 \end{array} \right.
\end{eqnarray}
with the boundary and initial conditions
 \begin{eqnarray}
 \nonumber  &v^i(0)=w(0)=w_x(0)=0, K_2 w_{xx}(L) = -k_1 \dot w_x(L)&\\
   \nonumber &\alpha^i_1h_i v^i_x(L)=-s_i \gamma^i \dot v^i (L),\quad i=1,3,&\\
\label{d-son-dcpld}  & K_1 \ddot w_x(L) -K_2 w_{xxx}(L) =k_2\dot w(L).&
%&(v^1, v^3,  w, \dot v^1, \dot v^3,\dot w)(x,0)=(v^1_0,  v^3_0, w_0, v^1_1, v^3_1,   w_1).&
 \end{eqnarray}
The operator $\mc{A}_{\phi}:\mc{H}\to\mc{H}$ is the coupling between the layers defined as the following
\begin{eqnarray}\label{opB}\mc {A}_{\phi}{\bf y}=
 \left( \begin{array}{c}
 0_{3\times 1} \\
 \mc  M^{-1} \left(H G_2 ~\phi^2_x\right)   \\
G_2\phi^2/(h_1 \rho_1)\\
-G_2 \phi^2/(h_3 \rho_3)\\
   \end{array} \right)
 \end{eqnarray}
 where ${\bf y}=(w,u^1, u^3, \tilde w, \tilde v^1, \tilde v^3)$ and  $\phi^2=\frac{1}{h_2}\left(-u^1+u^3 + H u_x\right).$
 Let $E_d(t)$ be natural energy corresponding to the system (\ref{dbas-dcpld})-(\ref{d-son-dcpld}), i.e. $\phi^2\equiv 0$ in (\ref{Energy-nat}).
  \begin{thm} Let ${\mc A}_d+B$ be the infinitesimal generator of the semigroup corresponding to the solutions of (\ref{dbas-dcpld})-(\ref{d-son}). Then the semigroup $\{e^{({\mc A}_d+B) t}\}_{t\ge 0}$ is exponentially
stable in $\mc H$.
 \end{thm}

 \textbf{Proof:}
 %Note that the equations in (\ref{dbas-dcpld}) are completely decoupled.
 The exponential stability of the semigroup $e^{({\mc A}_d+B)t}$ follows from the exponential stability of wave equations  \cite{O-Hansen4} and the Rayleigh beam equation \cite{B-Rao}.

\begin{Lem} \label{compact} The operator $A_{\phi}:\mc{H}\to \mc{H}$ defined in (\ref{opB}) is compact.
\end{Lem}
\vspace{0.1in}

 When $(w,u^1, u^3, \tilde w, \tilde u^1, \tilde u^3)\in \mathrm H,$ we have $w\in H^2_L(0,L)$ and $u^1, u^3 \in (H^1_L(0,L))^{2},$ and therefore $\phi^2 \in H^1_L(0,L).$ Since $\mc M: H^2_L(0,L)\to \mathrm L^2(0,L)$ remains  an isomorphism, the last terms in (\ref{opB}) satisfy
$ \mc M^{-1} \left(\phi^2_x\right)  \in H^2_L(0,L)$  where $\phi^2 \in H^1_L(0,L),$
and  $H^2_L(0,L)\times  H^1_L(0,L)$ is compactly embedded  in $H^1_L(0,L)\times (\mathrm L^2(0,L))^{2}.$ Hence the operator $A_{\phi}$  is compact in $\mc H$.

 \begin{thm}  Then the semigroup $\{e^{({\mc A}+B) t}\}_{t\ge 0}$ is exponentially
stable in $\mc H.$
 \end{thm}

 \textbf{Proof:} The semigroup $\mc A+B={\mc A}_d +B + {\mc A}_{\phi}$ is strongly stable on $\mc H$ by Theorem \ref{strong}, and the operator ${\mc A}_{\phi}$ is a compact in $\mc H$ by
Lemma \ref{compact}. Therefore, since the semigroup generated by $({\mc A}_d +B + {\mc A}_{\phi})-{\mc A}_{\phi}$ is uniformly
exponentially stable in $\mc H$ then the semigroup $\mc A=({\mc A}_d + B+ {\mc A}_{\phi})$ is uniformly
exponentially stable in ${\mc H}$ by e.g., the perturbation theorem of \cite{Trigg}.
\section{Future research}
%This paper uses a thorough variational approach to derive a model for an ACL beam with two piezoelectric layers by including magnetic effects. In contrast to the classical counterparts, our model is simplified substantially by ignoring the weight and the stiffness of the viscoelastic layer. We show that the system with magnetic effects is not strongly stable for certain choices material parameters if the feedback controllers for the longitudinal dynamics are chosen to be only electrical. As the magnetic effects are released, we  obtain the reduced Rao-Nakra type beam equations. We show that if we choose velocity type feedback controllers (for $V^{\rm T}$ and $V_{\rm B}$) for the longitudinal dynamics, and two additional feedback controllers for the transverse dynamics, the system is shown to be uniformly exponentially stable by using a compact perturbation argument.

% The Rayleigh beam equation $m \ddot w - K_1 \ddot{w}_{xx} + K_2 w_{xxxx} =0$ is still exponentially stable with only one feedback controller ($k_1\ne 0$, $k_2\equiv 0$ in (\ref{dbas-dcpld})), see \cite{B-Rao}. If we choose $B\equiv 0,$ i.e. $k_2\equiv 0,$ the exponential stabilizability of the model (\ref{dbas})-(\ref{d-son})  is currently unknown. The result all depends on whether the unique continuation result (Lemma \ref{xyz}) holds with only three boundary conditions for the bending at $x=L$  versus the classical four in this paper.

%Reducing the number of controllers to one for the bending motion of the whole composite is an open problem. This is currently under investigation by using spectral techniques.

A relevant  research problem  under investigation is whether we can recover the polynomial stability for certain combinations of material properties and for more regular initial data for the composite (\ref{dbas-eig-mag})-(\ref{d-son-eig-mag}).
For a single piezoelectric beam, this question is answered in \cite{Ozkan1}.% There are mathematical hurdles due to the complexity of coupling of the ACL composite.

The stabilization results obtained in this paper can be compared to the ones corresponding to the charge or current actuation (\cite{O-M2,Ozkan5}).

%There are ACL plate models with magnetic effects under consideration \cite{Ozkan2}. The exponential stabilizability of these models are more technical since the underlying decoupled dynamics are the Lam$\acute{e}$ system  coupled to the magnetic equations for the stretching motion, and the Kirchhoff's plate equation for the bending motion. This is currently under investigation.

\ifCLASSOPTIONcaptionsoff
  \newpage
\fi




%\vfill


\end{document}